\documentclass[letter]{amsart}

\usepackage{amssymb}

\newtheorem{thm}{Theorem}
\newtheorem{cor}[thm]{Corollary} 
\newtheorem{lem}[thm]{Lemma}
\newtheorem{prop}[thm]{Proposition}
\newtheorem{ex}[thm]{Example}

\newtheorem{rem}[thm]{Remark}

\title{$q$-Analogs of symmetric function operators}
\author{Michael Zabrocki}
\email{zabrocki@mathstat.yorku.ca}
\address{Mathematics and Statistics, York University,
Toronto, Ontario, M3J 1P3}
\address{ http://www.math.yorku.ca/\~{}zabrocki}

\begin{document}
\begin{abstract}
For any homomorphism $V$ on the space of symmetric functions,
we introduce an operation that creates a $q$-analog of $V$. By giving 
several examples we demonstrate that this quantization occurs naturally
within the theory of symmetric functions.  In particular, we show that the Hall-Littlewood
symmetric functions are formed by taking this $q$-analog
of the Schur symmetric functions
and the Macdonald symmetric functions appear by taking the $q$-analog of
the Hall-Littlewood
symmetric functions in the parameter $t$.  This relation is then used to derive 
recurrences on the Macdonald $q,t$-Kostka coefficients.

\noindent{\sc R\'esum\'e.} 
Pour un homomorphisme $V$ sur l'espace des fonctions sym\'etriques, 
nous pr\'esentons une op\'eration qui
cr\'ee un $q$-analogue de $V$. En donnant plusieurs exemples nous d\'emontrons que cette
quantization se produit naturellement dans la th\'eorie de fonctions sym\'etriques. En particulier,
nous prouvons que les fonctions sym\'etriques de Hall-Littlewood sont constitu\'ees en prenant ce
$q$-analogue des fonctions sym\'etriques de Schur et les fonctions sym\'etriques de Macdonald
apparaissent en prenant le $q$-analogue des fonctions sym\'etriques de 
Hall-Littlewood dans le param\`etre
$t$. Cette relation est alors employ\'ee pour d\'eriver des r\'ecurrence 
sur les coefficients Macdonald $q,t$-Kostka. 
\end{abstract}
\maketitle

\font\Sc=cmcsc10
\font\Ch=msbm9 scaled\magstep2
\font\gros=cmbx10 scaled\magstep2
\font\super=cmbx12 scaled\magstep2
\font\small=cmr8
\font\smallit=cmsl8
\font\smallb=cmbx8
\font\ninerm=cmr9

\newdimen\Squaresize \Squaresize=14pt
\newdimen\Thickness \Thickness=0.5pt

\def\Square#1{\hbox{\vrule width \Thickness
   \vbox to \Squaresize{\hrule height \Thickness\vss
      \hbox to \Squaresize{\hss#1\hss}
   \vss\hrule height\Thickness}
\unskip\vrule width \Thickness}
\kern-\Thickness}

\def\Vsquare#1{\vbox{\Square{$#1$}}\kern-\Thickness}
\def\Blk{\omit\hskip\Squaresize}

\def\Young#1{
\vbox{\smallskip\offinterlineskip
\halign{&\Vsquare{##}\cr #1}}}

\newdimen\squaresize \squaresize=5pt
\newdimen\thickness \thickness=0.2pt

\def\square#1{\hbox{\vrule width \thickness
   \vbox to \squaresize{\hrule height \thickness\vss
      \hbox to \squaresize{\hss#1\hss}
   \vss\hrule height\thickness}
\unskip\vrule width \thickness}
\kern-\thickness}

\def\vsquare#1{\vbox{\square{$#1$}}\kern-\thickness}
\def\blk{\omit\hskip\squaresize}
\def\noir{\vrule height\Squaresize width\Squaresize}%

\def\young#1{
\vbox{\smallskip\offinterlineskip
\halign{&\vsquare{##}\cr #1}}}

\def\thisbox#1{\kern-.09ex\fbox{#1}}
\def\downbox#1{\lower1.200em\hbox{#1}}

\def\la{{\lambda}}
\def\coeff{{\Big|}}
\def\scoeff{{\big|}}
\def\hs{\hskip .1in}
\def\hsk{\hskip .5in}
\def\endofproof {\hskip .1in $\diamondsuit$}
\def\keq{\simeq_k}
\def\snake {\rfloor}
\def\H {{\cal H}}
\def\V {{\cal V}}
\def\Sw {{\tilde S}}
\def\w{v}
\def\Part{\mathcal{P}}
\def\HHH {{\Ch H}}
\def\X {\hbox{\bf X}}
\def\Hw {{\tilde H}}
\def\sgn{{\varepsilon}}
\def\bproof{\noindent {\Sc Proof: }}
\newcommand{\ahat}[1]{\overline{#1}}
\newcommand{\htq}[1]{\widetilde{#1}^q}
\newcommand{\hqt}[1]{\widetilde{#1}^t}
\newcommand{\qht}[1]{{}^q\widetilde{#1}}
\newcommand{\DIW}[1]{\Z_\ge^{#1}}
\newcommand{\Pa}[1]{{\Ch P}^{#1}}
\newcommand{\Z}{{\mathbb Z}}

\def\eqnope{{}}

\def\pointir{\discretionary{.}{}{.\kern.35em---\kern.7em}\nobreak
\hskip 0em plus .3em minus .4em }
\def\article#1|#2|#3|#4|#5|#6|#7|
    {{\leftskip=7mm\noindent
     \hangindent=0mm\hangafter=1
     \llap{$[#1]$\hskip.35em}{#2,}
     #3, {\it #4} \nobreak {\bf #5} \nobreak {#6},
     \nobreak\ #7.\par}}
\def\unarticle#1|#2|#3|#4|#5|
    {{\leftskip=7mm\noindent
     \hangindent=0mm\hangafter=1
     \llap{$[#1]$\hskip.35em}{#2,}
     #3, {\it #4}\nobreak #5.\par}}
\def\livre#1|#2|#3|#4|
    {{\leftskip=7mm\noindent
    \hangindent=0mm\hangafter=1
    \llap{$[#1]$\hskip.35em}{#2,}
    {``#3,''} #4.\par}}

\section{ Introduction}

The Hall-Littlewood and Macdonald symmetric functions
are two examples of families of symmetric functions that
depend on a parameter $q$ such that setting this parameter
$q$ equal to $0$ yields one class of symmetric functions
which is not a product of generators
and setting the parameter $q$ equal to $1$ yields a
multiplicative basis.  There are other other classes of
symmetric functions with the same property, and in this article
we will show that practically any of these families are instances
of the same $q$-twisting of the symmetric function found by setting $q=0$.

This remarkable fact has lead to a completely elementary proof of
the polynomiality of the $q,t$-Kostka coefficients [GZ] and in this
article we use the very same observation to derive a combinatorial recurrence
on these coefficients as well as algebraic formulas for operators that add
a column to the partition indexing a Macdonald symmetric function.

The first section of this article will introduce some necessary notation
and the definition of this $q$-analog.  In the second section we give several
examples where it arises.  Some examples
will be nothing more than showing that $\ahat{V}$ for some $V$ 
is a formula that is well known in the literature.  Other examples
present some completely new equations, the most important of which will
concern the relation of the Hall-Littlewood symmetric functions to the
Macdonald symmetric functions.  This section will show that
this single $q$-analog appears in the creation of several different
classes of Schur positive symmetric functions.

In the third section we derive some formulas related to an operator that
adds a column to the Hall-Littlewood symmetric functions.  The $q$-analog
of this operator adds a column to the Macdonald
symmetric functions.  In the fourth section these equations are used
to give a formula for the action of this operator on the
Schur basis giving a combinatorial rule for computing the Macdonald
symmetric functions (i.e. a `Morris-like' recurrence for the $q,t$-Kostka
coefficients).

\section{Notation}

A partition of $n$ is a sequence of non-negative integers
$\la = (\la_1 \geq \la_2 \geq \la_3 \geq \cdots)$ such that
$\sum_i \la_i = n$.  The length of a partition is the largest index
$i$ such that $\la_i$ is nonzero, and it will be denoted here by $\ell(\la)$.
A partition will be drawn as a sequence of rows of boxes aligned at the
left edge with $\la_i$ cells
in the $i^{th}$ row.  We will use the French convention and draw these
diagrams with the largest row on the
bottom and the smallest row on the top. The conjugate partition $\la'$ is
the sequence whose $i^{th}$ entry is the number of cells in the $i^{th}$
column of the diagram for $\la$.

 The partition will be sometimes be 
identified with its diagram in the sort of language that is used.  For instance,
the operations of adding rows or columns to
partitions indexing bases for the symmetric functions are important here.  
the notation $(m,\mu)$ is used to represent the sequence $\mu$ with a part of
size $m$ prepended, which will be a partition as long as $\mu_1 \leq m$.
The notation $1^m|\mu$ will be used to represent the partition $(\mu_1+1,
\mu_2+1,\ldots, \mu_m+1)$ (as long as $\ell(\mu) \leq m$).

Let $\Lambda$ be the space of symmetric functions with the standard bases
for this space, $h_\la$ homogeneous, $e_\la$ elementary, $m_\la$
monomial, $f_\la$ forgotten, $p_\la$ power, and $s_\la$ the Schur
symmetric functions defined as they are in [M].  The involution $\omega$
that sends $p_k$ to $(-1)^{k-1} p_k$ relates these bases by $\omega(
h_\la ) = e_\la$, $\omega(m_\la) = f_\la$ and $\omega(s_\la) = s_{\la'}$.
The standard inner product on this space determines the 
dual bases
\begin{equation}
\left< p_\la, {p_\mu / z_\mu} \right> =
\left< h_\la, m_\mu \right> =
\left< e_\la, f_\mu \right> =
\left< s_\la, s_\mu \right> = \delta_{\la\mu},
\end{equation}
where $z_\la = \prod_{i>1} i^{n_i(\la)} n_i(\la)!$ 
with $n_i(\la)$ equal to the number
of parts of size $i$ in $\la$, and we have set
$\delta_{\la\la} = 1$ and $\delta_{\la\mu} = 0$ if $\la \neq \mu$.

For any element $f$ of $\Lambda$, let $f^\perp$ be the operation
that is dual to multiplication by $f$ with respect to the standard
inner product.  By definition we have that for any dual bases
$\{ a_\la \}_\la$ and $\{ b_\la \}_\la$, the action of $f^\perp$ on
another symmetric function $g$ is given by the formula
\begin{equation} f^\perp g = \sum_\la \left< g, f a_\la \right> b_\la.\end{equation}

`Plethystic' notation is a device for expressing the substitution of
the monomials of one expression in a symmetric function.
Assume that $E$ is a formal series in a
set of variables $x_1, x_2, ...$ with possible special parameters
$q$ and $t$.
For $k \geq 1$, set $p_k[E]$ to be $E$ with $x_i$ replaced by  
$x_i^k$ and $q$ and $t$ replaced by $q^k$ and $t^k$ respectively, that is
\begin{equation}p_k[E(x_1,x_2, \ldots,q,t)] = E(x_1^k, x_2^k, \ldots;q^k,t^k)\end{equation}
For an arbitrary symmetric function $P$, $P[E]$ will represent the
the formal series found by expanding $P$ in terms
of the power symmetric functions and then substituting $p_k[E]$
for $p_k$.
More precisely, if the power sum expansion of the symmetric
function $P$ is given by $P = \sum_\la c_\la p_\la$  then $P[E]$ is given by the 
formula
\begin{equation}P[E] = \sum_\la c_\la p_{\la_1}[E] p_{\la_2}[E] \cdots p_{\la_{\ell(\la)}}[E]. 
\end{equation}

The symmetric functions in the infinite set of 
variables $x_1, x_2, x_3, \ldots$ will be denoted by $\Lambda^X$.  $\Lambda$
and $\Lambda^X$ are isomorphic and in this exposition we will identify the
two spaces when it is convenient.  In plethystic notation, the isomorphism
that identifies the two spaces is given by $f \mapsto f[X]$ where
$X = x_1 + x_2 + x_3 + \cdots$ since under this map $p_k$ is sent to
$x_1^k + x_2^k + x_3^k+ \cdots$.  Also use the notation 
$X_n = x_1 + x_2 + \cdots + x_n$ to represent when
a symmetric function being evaluated a finite set of variables.  For sets of
variables using other letters we use a similar convention.

The symbol $\Omega = \sum_{n \geq 0} h_n$ will represent a
special generating function, and we will use
the plethystic notation for symmetric functions with this expression
as well with the following identities.
\begin{equation}
\Omega[X+Y] = \Omega[X]\Omega[Y]
\end{equation}
\begin{equation}
\Omega[X] = \prod_i \frac{1}{1-x_i}
\end{equation}
\begin{equation}
\Omega[XY] = \sum_{\la} s_\la[X] s_\la[Y]
\end{equation}

Operators that have the property that they add a row or a column to the
partition indexing a symmetric function will be known as `creation operators.'
The creation operators that will be used repeatedly are those that add a 
row to the Schur symmetric functions (due to Bernstein, see [Ze, p.~69], [M, p.~96])
and the Hall-Littlewood symmetric functions (due to Jing, see [J], [G] or 
[M, p.~238]).  In the third section, the operator introduced in [Za1]
that adds a column to the partition indexing a Hall-Littlewood symmetric
function will be developed further.

Define the following involution on the space $Hom(\Lambda, \Lambda)$ that 
is a useful tool
for deriving identities within the theory of symmetric functions. 
Let $V$ be an element of $Hom(\Lambda, \Lambda)$ and $P \in \Lambda$.  We define
the flip of $V$ by the formula

\begin{equation}\ahat{V} P[X] = V^Y P[X - Y] \coeff_{Y=X}.\end{equation}
It seems
to arise naturally when one considers the sorts of operators that will concern
us here (see [Za1] and [Za2]).  

The degree of a symmetric function $P \in \Lambda$ is the highest power
of $z$ in $P[zX]$ and will be denoted by $deg(P)$.  If $P[zX] = z^{deg(P)} P[X]$
then we will say that $P$ is of homogeneous degree.

Use this involution here to define a $q$-twisting of a symmetric function operator.  Let $V$ once
again be an element of $Hom(\Lambda, \Lambda)$ and let $F^q$ be defined by
$F^q P[X] = P[X (1-q)]$. Our $q$-analog is defined when
it acts of the symmetric function $P \in \Lambda$ by the formula

\begin{equation}\htq{V} P[X] = V^Y P[q X + (1-q)Y ] \coeff_{Y=X} =
\ahat{\ahat{V}F^q}P[X]. \label{qhat}\end{equation}

It is easily seen that this $q$-analog has the following fundamental propery.
\begin{rem}\label{rem1}
Let $V$ be an element of $Hom(\Lambda, \Lambda)$
and create the $q$-twisting of this operator from formula (\ref{qhat}), $\htq{V}$, and act this
new operator on a symmetric function $P[X]$ to create
an expression such as
\begin{equation}\htq{V} P[X]\end{equation}
This $q$-analog has the property that when $q=0$, the expression becomes
\begin{equation}V P[X]\end{equation}
and if $q=1$, then it reduces to the product
\begin{equation}V(1) P[X]\end{equation}
\end{rem}

This paper is concerned with generalizations of the standard bases, 
the Hall-Littlewood and
Macdonald symmetric functions, which depend on additional parameters $q$ and
$t$.  There are two important scalar products on the symmetric functions  related to
these bases.  They are defined by
their values on the power symmetric basis.
\begin{equation}\left< p_\la, p_\mu \right>_{t} = \delta_{\la\mu} z_\la \prod_{i=1}^{\ell(\la)} 
{1-t^{\la_i}} \end{equation}
\begin{equation}\left< p_\la, p_\mu \right>_{qt} = \delta_{\la\mu} z_\la \prod_{i=1}^{\ell(\la)} 
(1-q^{\la_i})(1-t^{\la_i}) \end{equation}

The Macdonald symmetric functions $H_\mu[X;q,t]$ are defined
by the following three conditions.
\vskip .2in

\noindent
1. $\left< H_\la[X;q,t], H_\mu[X;q,t] \right>_{qt} = 0 $ if $\la \neq \mu$.

\noindent
2.  $F^t H_\mu[X;q,t] = \sum_{\la \leq \mu} c_{\la\mu} m_\la[X]$ for suitable
coefficients $c_{\la\mu}$ and the sum is over all partitions $\la$ that are
smaller than $\mu$ in the standard dominance order.

\noindent
3. $\left<H_\mu[X], h_n[X] \right> = t^{n(\mu)}$ where $n(\mu) = \sum_{i} (i-1) \mu_i$.
\vskip .2in

The expansion of the $H_\mu[X;q,t]$ basis
in the Schur basis for the symmetric functions defines the coefficients
$K_{\la\mu}(q,t)$, that is
 $H_\mu[X;q,t] = \sum_{\la \vdash |\mu|} K_{\la\mu}(q,t) s_\la[X]$.

The Hall-Littlewood basis is defined similarly with respect to the $\left<,
\right>_{t}$ scalar product; simply stated $H_\mu[X;t] = H_\mu[X;0,t]$.

The symmetric functions $H_\mu[X;q,t]$ and $H_\mu[X;t]$ are the two families
of symmetric functions that will interest us the most here.  A symmetric
function with the property that when expressed in terms of the Schur basis
their coefficients are polynomials in $q$ and $t$ with non-negative coefficients
will be called {\it Schur positive}.  The Macdonald and Hall-Littlewood
functions are just two examples of families with this property.

\section{Examples}

\subsection{Schur symmetric functions I }

In [Ze] an operator attributed to Bernstein that adds a row to the Schur 
function
is given by $S_m = \sum_{i \geq 0} (-1)^i h_{m+i} e_i^\perp$.  
The formula has a very convenient form when expressed in terms
of plethystic notation.  Let $P[X]$ be a symmetric function in the $X$
variables.  Define a generating function of operators given by
\begin{equation}S(z) P[X] = P\left[X - {1 \over z} \right] \Omega[zX].\label{Szdef}\end{equation}

Now for any $m \in \Z$, set $S_m P[X] = S(z) P[X] \scoeff_{z^m}$.  If $m \geq 
\mu_1$, then it easily follows that $S_m s_\mu[X] = s_{(m,\mu)}[X]$.
$S_m$ is a creation 
operator for the Schur basis since we have the
formula
\begin{equation}S_{\mu_1} S_{\mu_2} \cdots S_{\mu_{\ell(\mu)}} 1 = s_{\mu}[X].\end{equation}

Now if we set $H(z) = \htq{S(z)}$ and $H_m^q = H(z) \scoeff_{z^m}$, 
then this is a $q$-analog of the operator
$S_m$ and we may calculate that
\begin{eqnarray}
H(z) P[X] &=& S(z) P[q X + (1-q)Y ] \coeff_{Y=X}\nonumber\\
&=& P\left[q X + (1-q)\left(X - {1 \over z}\right) \right] \Omega[zX] \label{StoH}\\
&=&
 P\left[ X - {1-q \over z} \right] \Omega[zX]\nonumber
\end{eqnarray}

Remarkably, this is the formula for the Hall-Littlewood creation operator
of Jing [J] in the notation used by Garsia [G].  These operators
have the property that

\begin{thm} (Jing [J])  Let $H_m^q = \htq{S_m}$. Then
\begin{equation}H_{\mu_1}^q H_{\mu_1}^q \cdots H_{\mu_{\ell(\mu)}}^q 1 = H_{\mu}[X;q].
\end{equation}
\end{thm}

The fact that $H_\mu[X;0] = s_\mu[X]$ and $H_\mu[X;1] = h_\mu[X]$ follows from
Remark \ref{rem1}.

\subsection{ Schur symmetric functions II }\label{Schur2}

For any sequence of integers $\gamma = (\gamma_1, \gamma_2, \ldots, \gamma_n)$,
let $x^\gamma$ represent the monomial $x_1^{\gamma_1} x_2^{\gamma_2} \cdots
x_n^{\gamma_n}$.
We say that $\gamma$ is a dominant weight if $\gamma_1 \geq \gamma_2 \geq \cdots
 \geq \gamma_n$.

  The symmetric group $S_n$ acts on any polynomial in the $X_n$
variables by permuting their indices.  For any $\sigma \in S_n$, set $\sgn_\sigma$
to be the sign of the permutation.  Also set $\delta = (n-1, n-2, \ldots, 1,0)$.
Then define for any polynomial $f$ in the $x_i$ variables a symmetrization
operator $\pi_n (f) = J(f) / J(1)$ where

\begin{equation}J(f) = \sum_{\sigma \in S_n} \sgn_\sigma \sigma(x^\delta f).\end{equation}
When $\la$ is a partition, $\pi_n$ sends $x^\lambda$ to the Schur function
$s_\la[X_n]$.  When
$\la$ is any dominant weight, set $s_\la(X_n) = \pi_n( x^\la)$ to be
the resulting Laurent polynomial.

Let $\eta$ be a sequence of positive integers whose sum is $n$ (a composition)
and define the set of ordered pairs 
$Roots_\eta = \{ (i,j) : 1 \leq i \leq \eta_1+\eta_2+ \cdots+\eta_r < j \leq n
\hbox{ for some }r\}$.  Consider the following formal power series given
by the formula

\begin{equation} H_{\mu\eta}(X_n; q) = \pi_n \left(x^\mu\prod_{(i,j) \in Roots_\eta} 
(1-q x_i/x_j)^{-1}\right).\end{equation}
This formal series has an expansion in terms of the Schur functions indexed by
all dominant integral weights.  Define $K_{\la\mu\eta}(q)$ as the coefficient
of $s_\la(X_n)$ in $H_{\mu\eta}(X_n; q) $, so that
\begin{equation}H_{\mu\eta}(X_n; q) = \sum_{\la} K_{\la\mu\eta}(q)s_\la(X_n),\end{equation}
where the sum is over all dominant integral weights $\la$.  The
$K_{\la\mu\eta}(q)$ are known as the generalized or parabolic Kostka polynomials.
For a more complete exposition of these polynomials we refer the reader to
$[KS]$, $[SW]$, or $[K]$.

Now consider a composition of the Bernstein Schur function operators $S_\nu =
S_{\nu_1} S_{\nu_2} \cdots S_{\nu_{\ell(\nu)}}$ and define $H_{\nu}^q = \htq{S_\nu}$.
By a calculation similar to (\ref{StoH}), one may show that

\begin{equation}H_\nu^q P[X] = P \left[ X - (1-q) Z^* \right] \Omega[ZX] \prod_{1 \leq j < i \leq \ell(\nu)} 1-z_i/z_j
\coeff_{z_1^{\nu_1} z_2^{\nu_2} \cdots z_{\ell(\nu)}^{\nu_{\ell(\nu)}}}, \end{equation}
where $Z^* = \sum_{i=1}^{\ell(\nu)} 1/z_i$.

In work with Mark Shimozono [SZ], we demonstrated that $H_\nu^q$ is an operator with
interesting properties related to the 
generalized Kostka coefficients.  In particular, they can be used as generating
functions for the generalized Kostka coefficients.  They are also operators which
act on symmetric functions and can be used to build a family of symmetric functions.

\begin{thm} (Shimozono, Zabrocki [SZ])  Let $\eta$ be a composition of
$k$ and $\mu \in \Z^k$,
set $\mu^{(i)} = (\mu_{\eta_1+\cdots + \eta_{i-1}}, \ldots,
\mu_{\eta_1+\cdots + \eta_{i}})$.
For any $\nu \in \Z^\ell$ set $H_\nu = \htq{S_\nu}$ where
$S_\nu = S_{\nu_1} \cdots S_{\nu_{\ell(\nu)}}$, then we have
\begin{equation}H_{\mu^{(1)}}^q H_{\mu^{(2)}}^q  \cdots H_{\mu^{(\ell)}}^q =
\sum_{\la}  K_{\la\mu\eta}(q) H_\la^q,\end{equation}
where the sum is over all dominant weights $\la$.
In particular, when this operator is applied to the symmetric function $1$,
we arrive at  a class of symmetric functions and may set
\begin{equation} H_{\mu\eta}[X;q] :=
H_{\mu^{(1)}}^q H_{\mu^{(2)}}^q  \cdots H_{\mu^{(k)}}^q 1 =
\sum_{\la}  K_{\la\mu\eta}(q) s_\la[X],\end{equation}
where the sum is over all partitions $\la$.
\end{thm}

In addition, these operators also seem to be fundamentally related to the
new class of symmetric functions referred to as `atoms' $A_\la^{(k)}[X;q]$
(see [LLM], [LM1], [LM2]).  In particular,
when the partition indexing the operator is a rectangle, it was conjectured
that $H_{(\ell^{k+1-\ell})}^q$ is a creation operator for this class of
symmetric functions.

By Remark 1, the  coefficients 
$\left< H_{\mu\eta}[X;q], s_\la[X] \right>= K_{\la\mu\eta}(q)$
have the property that when $q=1$,
they are the Littlewood-Richardson coefficients 
$c_{\mu^{(1)}\mu^{(2)}\cdots\mu^{(k)}}^\la$ and when $q=0$, the 
$K_{\la\mu\eta}(0)=1$ if $\mu=\la$ and $0$ otherwise.  It is conjectured that the coefficients
$K_{\la\mu\eta}(q)$ are polynomials in $q$ with non-negative integer
coefficients.

\subsection{ Homogeneous creation operator }

Multiplication by $h_k$ is an operator that adds a row to the homogeneous
symmetric functions.  The $q$-analog of this operator is once again $h_k$
and so this is not particularly interesting.  However, in [Za2] we gave a
formula for an operator that adds a column to the homogeneous symmetric 
functions and in [Za1] we gave a combinatorial description of the action of
this operator on the Schur function basis.  Let $H_{1^m}$ be a family of
operators with the property that
\begin{equation}H_{1^{\la_1}} H_{1^{\la_2}} \cdots H_{1^{\la_{\ell(\la)}}} 1 =
h_{\la'}[X].\end{equation}

The $q$-twisting of this operator is another example where this $q$-analog appears
naturally to produce a family
of Schur postitive symmetric functions.  It develops that
$\htq{H_{1^m}}$ is an operator that adds a column to the
the symmetric functions $(q;q)_\la h_{\la}\left[{X \over 1-q}\right]$, where
$(q;q)_k = (1-q)(1-q^2)\cdots (1-q^k)$ and $(q;q)_\la = (q;q)_{\la_1} 
(q;q)_{\la_2}\cdots (q;q)_{\la_{\ell(\la)}}$.  This family has the property
that when $q=1$ the symmetric functions become $h_{1^{|\la|}}$ and when
$q=0$ they become $h_\la$.

\begin{thm}   Let $H_{1^m}$ be an operator
with the property that $H_{1^m} h_\la[X] = h_{1^m|\la}[X]$ for 
$\ell = \ell(\la) \leq m$.  Then we have
\begin{equation}\htq{H_{1^{\la_1}}} \htq{H_{1^{\la_2}}} \cdots \htq{H_{1^{\la_\ell}}} 1 = 
(q;q)_{\la'} h_{\la'} \left[ {X \over 1-q} \right]\end{equation}
\end{thm}


\bproof
\begin{equation}\htq{H_{1^m}}\left( (q;q)_\la h_{\la}\left[{X \over 1-q}\right]\right) =
(q;q)_\la H_{1^m}^Y \left(h_{\la}\left[{qX + (1-q)Y \over 1-q}\right]\right) 
\coeff_{Y=X }\label{qhathomo}\end{equation} 

Since the summation formula for a homogeneous symmetric function in two sets of
variables is given as $h_m[X+Y] = \sum_{i=0}^m h_i[X] h_{m-i}[Y]$, then for a
partition $\la$,
$h_\la[X+Y] = \prod_{n=1}^{\ell(\la)}\sum_{i=0}^{\la_n} h_{i}[X] h_{\la_n-i}[Y]$. 
 It follows
that (\ref{qhathomo}) reduces to 
\begin{eqnarray}
&=&
(q;q)_\la \prod_{n=1}^{m}\sum_{i=0}^{\la_n} 
 h_{i}\left[{qX \over 1-q}\right] h_{\la_n-i+1}[X]\nonumber\\
&=& (q;q)_\la \prod_{n=1}^{m}\left(\sum_{i=0}^{\la_n+1} 
 h_{i}\left[{qX \over 1-q}\right] h_{\la_n-i+1}[X] - 
h_{\la_n+1}\left[{qX \over 1-q}\right]\right)\nonumber\\
&=& (q;q)_\la \prod_{n=1}^{m} h_{\la_n+1}\left[{X \over 1-q}\right]
 - q^{\la_n+1} h_{\la_n+1}\left[{X \over 1-q}\right]\\
&=& (q;q)_\la \prod_{n=1}^{m} (1-q^{\la_n+1}) h_{\la_n+1}\left[{X \over 1-q}\right]\nonumber\\
&=& (q;q)_{1^m|\la} h_{1^m|\la}\left[{X \over 1-q}\right].\nonumber\end{eqnarray}
\endofproof

There are elementary proofs that the functions
$(q;q)_\la h_{\la}\left[{X \over 1-q}\right]$ are Schur positive.  
Once again we have a case of a Schur positive 
$q$-analog arising from
a Schur positive family of symmetric functions
$h_\la[X]$.
Remark \ref{rem1} implies that the limit as $q$ goes to $1$ of these
symmetric functions is $h_{1^n}[X]$ when $\la$ is a partition of $n$
and when $q=0$ we have that 
they reduce to $h_{\la'}[X]$.

\subsection{ Macdonald's Operators }

Macdonald introduced operators $D_n^r$ (see [M] p. 315) such that the Macdonald
polynomials $P_\la(X;q,t)$ are characterized as eigenfunctions of this family of 
operators.
What we will show is that the Macdonald operators are the $q$-twisting of the
same operators with $q$ set equal to $0$.

Let $T_{x_i}$ be an operator on polynomials with the property $T_{x_i} P[X_n] =
P[X - x_i]$.  Consider the operator
\begin{equation}D_n^r(t) = \sum_{I \subseteq \{1, \ldots, n\}}A_I(X_n;t) \prod_{i \in I} T_{x_i}\end{equation}
where $A_I(X_n;t) = t^{r \choose 2} \prod_{i \in I, j \notin I} {t x_i - x_j \over 
x_i -x_j}$ and the sum is over all subsets $I$ of size $r$.

Now define $D_n^r(q,t) = \htq{ D_n^r(t)}$, then calculate that
\begin{eqnarray}
D_n^r(q,t) P[X_n] &=& 
\sum_{I \subseteq \{1, \ldots, n\} \atop |I| = r}A_I(Y_n;t)\prod_{i \in I} T_{y_i} 
P[q X_n + (1-q)Y_n] \coeff_{Y_n = X_n}\nonumber\\
&=& \sum_{I \subseteq \{1, \ldots, n\} \atop |I| = r}A_I(Y_n;t)  P[X_n - (1-q)X_I].
\end{eqnarray}
If we set $T_{q,x_i} P[X_n] = P[X_n - (1-q)x_i]= \htq{T_{x_i}}$, then
\begin{equation}D_n^r(q,t) = \sum_{I \subseteq \{1, \ldots, n\} \atop |I|=r}A_I(X_n;t) 
\prod_{i \in I} T_{q,x_i},\label{macop}\end{equation}
and these are the operators $D_n^r$ as they are defined in [M].  As was presented
in this reference, set 
$D_n(u;q,t) = \sum_{r=0}^n u^r D_n^r(q,t)$.

Consider again the Bernstein operators, $S_m$, as they were defined in
equation (\ref{Szdef}) and also consider
${\tilde S}_m = \omega S_m \omega = (-1)^m S(-z) \scoeff_{z^m}$.  From these
operators create a $q,t$ analog by applying this parameter deformation
twice.  Set ${\tilde D}_m = {\widetilde{\htq{{\tilde{S}}_m}}}^{1/t}$ and ${\tilde D}_m^* 
= {\widetilde{{\widetilde{S_m}}^{t}}}^{1/q}$.

A simple calculation yields that
\begin{equation}{\tilde D}_m P[X] = P\left[ X + {(1-q)(1-1/t) \over z}\right] 
\Omega[-zX] \coeff_{z^m}, \end{equation} and
\begin{equation}{\tilde D}_m^* P[X] = P\left[ X - {(1-1/q)(1-t) \over z}\right] 
\Omega[zX] \coeff_{z^m}. \end{equation}

These families of operators were studied in [GHT] and more extensively in 
[BGHT] to show the polynomiality of the $q,t$-Catalan numbers.  In particular, when
$m=0$ it is known that these operators are related to the operators $D_n^1$ 
and that they have the family $H_\mu[X;q,t]$ as eigenfunctions.  Theorem $1.2$
of [GHT] was the following result (translated on the $H_\mu[X;q,t]$ basis).

\begin{thm} (Garsia-Haiman-Tesler [GHT]) For $\mu$ a partition of $n$, we have
\begin{equation}{\tilde D}_0 H_\mu[X;q,t] = \Bigg(1 - (1-1/t) \prod_{i \geq 1} t^{-i} 
(1-q^{\mu_i})\Bigg)
H_\mu[X;q,t],\end{equation}
\begin{equation}{\tilde D}_0^* H_\mu[X;q,t] = \Bigg(1 - (1-t) \prod_{i \geq 1} t^{i} 
(1-q^{-\mu_i})\Bigg)
H_\mu[X;q,t].\end{equation}
\end{thm}

\subsection{ Hall-Littlewood creation operator }

Consider now an operator that adds a column to the Hall-Littlewood symmetric
functions $H_\la[X; t]$.  One such operator was 
introduced in [Za2] where the
combinatorial action on the Schur function basis was discussed and some explicit
formulas were presented.  We will discuss some of these operators in more detail in 
the following section.  For now we present the following theorem.

\begin{thm} 
\label{Hqtcreation}
Let $\la$ be a partition such that $ \ell = \ell(\la) \leq m$.
Any operator ${H^{t}_{1^m}}$ with the property
${H^{t}_{1^m}} H_\mu[X;t] = H_{1^m|\la}[X;t]$ satisfies the equation
\begin{equation}\htq{H^{t}_{1^{\la_1}}} \htq{H^{t}_{1^{\la_2}}} \cdots 
\htq{H^{t}_{1^{\la_\ell}}}  = H_{\la'}[X;q,t]\end{equation}
\end{thm}

The most elementary proof of this theorem can be seen in [GZ] and is
self contained and uses nothing
more than the identities in the original paper of Macdonald.
We present here a short proof that
follows by demonstrating that if the theorem is true for one such
operator $H_{1^m}^t$, then it is true for all such operators. 

 Some of the first proofs of the polynomiality
of the $q,t$-Kostka coefficients used operators of this type to show that
the $q,t$-Kostka polynomials satisfied recurrences that did not have
denominators.  Almost any of the operators
given by Kirillov-Lapointe-Noumi-Vinet [KN1] [KN2]  [LV1] [LV2] 
are of this type.  We need only one
example and the following lemma.

\begin{lem} \label{allwork} Let $\la$ be a partition such that $\ell(\la) \leq m$.
Assume that there exists some operator $H_{1^m}^{t}$ such that 
\begin{equation}H_{1^m}^{t} H_\la[X;t] = H_{1^m|\la}[X;t] \label{eqa}\end{equation}
 and 
\begin{equation}\htq{H_{1^m}^t} H_\la[X;q,t] = H_{1^m|\la}[X;q,t]. \label{eqb}\end{equation}
Then for any operator ${H'}^{t}_{1^m}$ that satisfies equation (\ref{eqa})
will also satisfy equation (\ref{eqb}).
\end{lem}

\bproof It follows from the definition of the Macdonald symmetric functions
and the property that $H_\mu[X;t,q] = \omega H_{\mu'}[X;q,t]$
that we have the triangularity relation
$H_\mu[X(1-q);q,t]= \sum_{\la \geq \mu}
a_{\la\mu}(q,t) s_\la[X]$. Since $H_\mu[X;t]= \sum_{\la \geq \mu}
K_{\la\mu}(t) s_\la[X]$, then there exist coefficients $b_{\la\mu}(q,t)$
such that $H_\mu[(1-q)X;q,t] = \sum_{\la \geq \mu} b_{\la\mu}(q,t) H_\la[X;t]$.

Consider the expansion of $H_\mu[X+Y;q,t] = \sum_{\nu \subseteq \mu}
H_{\mu \slash \nu}[X;q,t] H_\nu[Y;q,t] $ which may be seen as a transformation
of formula $(7.9)$ on p.~345 of [M].

Since our operator $H_{1^m}^t$ adds a column of size $m$ to $H_\nu[X;t]$ and 
$\htq{H_{1^m}^t}$ adds a column to $H_\mu[X;q,t]$, then 
$H_{1^m|\mu}[X;q,t]$ is given by the formula

\begin{eqnarray}\htq{H_{1^m}^t} H_\mu[X;q,t] &=&
\sum_{\la \subseteq \mu} H_{\mu \slash \la}[q X;q,t] \;
H_{1^m}^t \, H_{\la}[(1-q) X;q,t]\nonumber \\
&=& 
\sum_{\la \subseteq \mu} H_{\mu \slash \la}[q X;q,t]
\sum_{\nu \geq \la} b_{\la\nu}(q,t)\; H_{1^m}^t \, H_\nu[X;t] \\
&=&
\sum_{\la \subseteq \mu} H_{\mu \slash \la}[q X;q,t]
\sum_{\nu \geq \la} b_{\la\nu}(q,t)\, H_{1^m|\nu}[X;t].\nonumber
\end{eqnarray}

The right hand side of this expression is independent of the operator that
is used to derive it, therefore absolutely any operator ${H'}^t_{1^m}$ that
has the property that ${H'}^t_{1^m} H_\mu[X;t] = H_{1^m|\mu}[X;t]$, also
has the property that $\htq{{H'}^t_{1^m}}$ adds a column to the Macdonald
symmetric functions $H_\mu[X;q,t]$. \endofproof

To prove the theorem it is necessary to produce at least one
operator that satisfies the property of the lemma.  Fortunately this is relatively
easy since practically any in the development in
[KN1] [KN2] [LV1] [LV2] such that one can set $q=0$ without needing to
take a limit satisfies the properties of Lemma \ref{allwork} (once transformed
to the $H_\mu[X;q,t]$ basis).

Consider Expression $3$ of [LV1].  Let $D_I^r(u;q,t)$ be the Macdonald
operator of equation (\ref{macop}) that acts only on the variables $x_i$ for
$i$ in the set $I$. A formula for an operator that adds a column of size
$m$ onto the $J_\mu[X;q,t] := H_\mu[X(1-t);q,t]$ basis is given by
\begin{equation}B_m^{(3)}(q,t) = \sum_{|I| = m} \sum_{r = 0}^m t^{-r} x_I \prod_{i \in I\atop
j\notin I} {x_i - x_j/t \over x_i - x_j} D_I^r(-t;q,t)\end{equation}
where $x_I$ here represents $\prod_{i \in I} x_i$.  It is easy to demonstrate
by acting on an arbitrary symmetric function that
\begin{equation}B_m^{(3)}(q,t) = \htq{ B_m^{(3)}(0,t)} \end{equation}

Let $F^t$ be an operator that sends the symmetric function $H_\mu[X;q,t]$
to the symmetric function $J_\mu[X;q,t]$.  More precisely, for an arbitrary
symmetric function $P[X]$  set $F^t P[X] = P[X (1-t)]$ and denote the
inverse of this operator $F_t^{-1}$.  Since $B_m^{(3)}(q,t)$
adds a column of size $m$ to the $J_\mu[X;q,t]$ basis, 
$F_t^{-1} B_m^{(3)}(q,t) F^t$
is an operator that adds a column to the $H_\mu[X;q,t]$ basis.

When $q=0$ in the operator $B_m^{(3)}(q,t)$ it becomes an operator that adds
a column to the $H_\mu[X(1-t);t]$ basis.  That is, we have
\begin{equation}B_m^{(3)}(0,t) H_\mu[X(1-t);t] = H_{1^m|\mu}[X(1-t);t].\end{equation}
Since $F^t H_\mu[X;t] = H_\mu[X(1-t);t]$, this implies $F_t^{-1} B_m^{(3)}(0,t) F^t$
is an operator that adds a column to the $H_\mu[X;t]$ basis.

To demonstrate the theorem, it remains to show that the $q$-twist of
$F_t^{-1} B_m^{(3)}(0,t) F^t$ is exactly the operator 
$F_t^{-1} B_m^{(3)}(q,t) F^t$.  This follows from the fact that conjugation
by $F^t$ commutes with the $q$-twisting for any symmetric function operator.

\begin{lem} For $V \in Hom(\Lambda, \Lambda)$ we have
\begin{equation}F_t^{-1} \htq{V} F^t = \htq{F_t^{-1} V F^t}.\end{equation}\end{lem}

\bproof This follows by acting both the left and the right hand side
of this equation on an arbitrary symmetric function. \endofproof

This leads us to several other formulas for operators with similar properties.
Consider the following corollary.

\begin{cor}  Define $h_\la(q) = \prod_{s \in \la} 1-q^{a_\la(s) + l_\la(s)+1}$.
Let $H^q_{1^m}$ be an operator with the property $H^q_{1^m} H_{\la}[X;q]=
H_{1^m|\la}[X;q]$ for $\ell = \ell(\la) \leq m$. Then
\begin{equation}\hqt{H^q_{1^{\la_1}}}  \hqt{H^q_{1^{\la_2}}} \cdots \hqt{H^q_{1^{\la_\ell}}} 1
= \omega H_{\la}[X;q,t], \end{equation}
and
\begin{equation}\htq{H^q_{1^{\la_1}}}  \htq{H^q_{1^{\la_2}}} \cdots \htq{H^q_{1^{\la_\ell}}} =
h_{\la'}(q) s_{\la'} \left[ {X \over 1-q} \right].\end{equation}
\end{cor}

\bproof This follows from Theorem \ref{Hqtcreation} and the following two identities
about Macdonald's symmetric functions.

\begin{equation}H_\mu[X;t,q] = \omega H_{\mu'}[X;q,t]\end{equation}
\begin{equation}H_\mu[X;q,q] = h_\mu(q) s_\mu\left[ {X \over 1-q} \right]\end{equation}
\endofproof

In the next two sections we will develop this last example in more detail.  The
observation that the Macdonald polynomials are built up from the $q$-twisting
of the operators that build the Hall-Littlewood symmetric functions allows
us to derive several interesting formulas for these operators.

\begin{rem}\label{rem2} Unfortunately, the
analog of section \ref{Schur2}
does not seem to extend to these operators as a way of generalizing
the Maconald symmetric functions.  Consider the
operator $H_{1^\la}^t := H_{1^{\la_1}}^t H_{1^{\la_2}}^t \ldots H_{1^{\la_{\ell(\la)}}}^t$. 
We would hope that a composition of $\htq{H_{1^\la}^t}$ are Schur
positive if reasonable conditions are placed on $\la$.
By calculating examples we begin to be encouraged by such
a conjecture, however for large enough examples
it seems to break down (for example if
$\la^{(1)} = (4), \la^{(2)} = (2,2), \la^{(3)} = (1,1)$, then
$\htq{H_{1^{\la^{(1)}}}^t} \htq{H_{1^{\la^{(2)}}}^t} \htq{H_{1^{\la^{(3)}}}^t} 1$
is not Schur positive).
\end{rem}

\section{Ribbons and Hall-Littlewood symmetric functions}

In [Za1] we gave a combinatorial formula for the action of an operator
that adds a column to the Hall-Littlewood symmetric functions.  We will
recall some of the definitions and theorems from that work and use them
to derive some useful formulas.

The definition of a ribbon is a skew partition that
contains no $2 \times 2$ subdiagrams.  For a non-empty partition $\la$, define
$\la^{rc} = (\la_2-1, \la_3-1, \ldots, \la_{\ell(\la)}-1)$ (the $rc$ 
indicates that $\la$ has the first row and first column removed).
If $R$ is a ribbon of size
$m$ (denoted by $R \models m$) then $R$ will be equal to $\la \slash \la^{rc}$ 
for some partition
$\la$ with $\la_1+\ell(\la) -1 = m$. 

Set $D(R)$ equal to the descent set of $R$, that is the set 
$\{ i~|~i+1^{st} $ cell lies below the $i^{th}$ cell in R $\}$
when the cells are labeled with the integers $1$ to $n$ from left to
right and top to bottom.  Therefore every ribbon can be identified
with a subset of $\{1, \ldots, m-1\}$.

There is a natural statistic associated with a ribbon.  Define the
major index of a ribbon to be 
$maj(R) = \sum_{i \in D(R)} i$.  Its complementary statistic will be
$comaj(R) = {|R| \choose 2} - maj(R)$.

From formula (\ref{Szdef}), $S_m$ is an
operator that adds a row to the partition indexing a Schur symmetric
function.  By conjugating $S_m$ by $\omega$, one obtains an operator
that adds a column.  Define $\Sw_m = \omega S_m \omega$.  In plethystic
notation, this operator is given as
\begin{equation}\Sw_m P[X] = (-1)^m P\left[ X + {1\over z}\right] \Omega[-zX] 
\coeff_{z^m}.\end{equation}

Now for each ribbon of size $m$, define an operator that raises the
degree of a symmetric function by $m$.  For $R = \la \slash \la^{rc}$ set
\begin{equation}S^R = s_{\la^{rc}}^\perp \Sw_{\la_1'} \Sw_{\la_2'} \cdots \Sw_{\la_{\la_1}'},\end{equation}
where $\la_i'$ is the length of the $i^{th}$ column in $\la$.
This is a combinatorial operator in the sense that all calculations can be
computed on the Schur basis 
using the Littlewood-Richardson rule and the commutation relations
$\Sw_a \Sw_b = - \Sw_{b-1} \Sw_{a+1}$ and $\Sw_a \Sw_{a+1} = 0$ so that the
operator $S^R$ can be thought of as an operator that acts on $s_\la$ by adding
the ribbon $R$ to the left of $\la$.

The main theorem in $[Za1]$ was the following result.

\begin{thm}{(Theorem 1.1 of [Za1])}\label{ribop} The operator $H_{1^m}^q = \sum_{R \models m} q^{comaj(R)} S^R$ has the
property that $H_{1^m}^q H_\mu[X;q] = H_{1^m|\mu}[X;q]$ for $\ell(\mu) \leq m$.
\end{thm}

Some elegant relations develop with the flip operation and ribbon operators.
Note that it follows directly from the definition that if $R$ is a ribbon
of size $m$ and $R^+$ is a ribbon of size $m+1$ with $D(R) = D(R^+)$, then 
$S^{R^+} = S^R \Sw_1$.
It develops that there is also a recursive method for adding a cell 
below the ribbon.  If $R_+$ is a ribbon of size $m+1$ such that 
$D(R_+) = D(R) \cup \{m\}$,
then we have the following surprising formula.

\begin{thm}{(Theorem 2.2 of [Za1])} \label{addacell} If $R \models m$ and $R_+ \models m+1$
such that $D(R_+) = D(R) \cup \{m\}$, then $S^{R_+} = \ahat{\ahat{S^R} S_1}$.
\end{thm}

This theorem can be used to produce the following plethystic formula for a ribbon operator.

\begin{prop}\label{plethrib} Let $R \models m$, then 
\begin{equation}S^R P[X] = (-1)^{m-|D(R)|} P[ X+ Z^*] \Omega[-(z_1 + Z_{D(R)^c+1}) X] \prod_{1 \leq i < j \leq m}
(1-z_j/z_i) \coeff_{z_1 z_2 \cdots z_m},\end{equation}
where we have set $Z^* = \sum_{i=1}^m 1/z_i$ and $Z_{D(R)^c+1} = 
\sum_{i \in [1,m-1]-D(R)} z_{i+1}$.
\end{prop}

\bproof By induction using Theorem \ref{addacell} and direct calculation. \endofproof

It follows from Theorem \ref{ribop} and \ref{addacell} that $H_{1^m}^t$
may be defined recursively.
Set $H_{1^1}^t = \Sw_1$ and 
\begin{equation}H_{1^{m+1}}^t = t^m H_{1^m}^t \Sw_1 + \ahat{\ahat{H_{1^m}^t} S_1}.\end{equation}
Either from this recursive definition or from the previous proposition,
one may demonstrate the following plethystic formula for the $H_{1^m}^t$ operator.

\begin{prop}  The operator $H_{1^m}^t$ of Theorem \ref{ribop}
has the following form in
plethystic notation.
\begin{equation}H_{1^m}^t P[X] = -P[ X+ Z^*]\Omega[-z_1 X] \prod_{i=2}^m (1 - t^{i-1} \Omega[-z_i X])
\prod_{1 \leq i < j \leq m}
(1-z_j/z_i) \coeff_{z_1 z_2 \cdots z_m} \end{equation}
\end{prop}

Note that because the coefficient of $z_1$ in the expression
\begin{equation}P[X+Z^*]\prod_{i=2}^m (1 - t^{i-1} \Omega[-z_i X]) \prod_{1 \leq i < j \leq m}(1-z_j/z_i)\end{equation}
is zero, we also have the following equivalent expression.

\begin{cor}\label{plethH1mt} The operator $H_{1^m}^t$ of Theorem \ref{ribop} has the following form in
plethystic notation.
\begin{equation}H_{1^m}^t P[X] = P[ X+ Z^*] \prod_{i=1}^m (1 - t^{i-1} \Omega[-z_i X])
\prod_{1 \leq i < j \leq m}
(1-z_j/z_i) \coeff_{z_1 z_2 \cdots z_m} \label{H1mt}\end{equation}
\end{cor}

In the next section Theorem \ref{Hqtcreation} will be used to develop methods for
computing Macdonald polynomials and the $q,t$-Kostka coefficients from
these formulas.
One may use some of the properties of the ribbon operators to derive
several other formulas for operators $H_{1^m}^t$ and hence for
$\htq{H_{1^m}^t}$, but this particular formula seems like a natural
extension to the ribbon operator formula for the Hall-Littlewood symmetric
functions.

\section{Generalized ribbons and Macdonald symmetric functions}

Consider the following generalization of the plethystic formulas
presented in the previous section.  
Since we know from Theorem \ref{Hqtcreation} that the operator
$\htq{H_{1^m}^t}$ is an operator that adds a column to the Macdonald
symmetric functions, the $q$-analog of equation (\ref{H1mt}) yields the following theorem.

\begin{thm}\label{plethH1mqt}  The following operator adds a column to the Macdonald symmetric
functions $H_{\mu}[X;q,t]$ if $ \ell(\mu) \leq m$.
\begin{equation}H_{1^m}^{qt} P[X] = P[ X+ (1-q)Z^*] \prod_{i=1}^m (1 - t^{i-1} \Omega[-z_i X])
\prod_{1 \leq i < j \leq m}
(1-z_j/z_i) \coeff_{z_1 z_2 \cdots z_m} \label{H1mqt}\end{equation}
\end{thm}

\bproof  This follows from Theorem \ref{Hqtcreation} and Corollary \ref{plethH1mt}.  
Calculate (\ref{H1mqt}) by using equation (\ref{H1mt}) and (\ref{qhat}) to
show
\begin{eqnarray}\htq{H_{1^m}^{t}} P[X] &=& H_{1^m}^{tY} P[q X + (1-q)Y] \coeff_{Y=X}\nonumber\\
&=& P[ q X + (1-q) Y+ (1-q) Z^*] \prod_{i=1}^m (1 - t^{i-1} \Omega[-z_i Y])\\
& &\hskip .5in \prod_{1 \leq i < j \leq m}
(1-z_j/z_i) \coeff_{z_1 z_2 \cdots z_m}\coeff_{Y=X} .\nonumber\end{eqnarray}
\endofproof

We will
develop this operator further and show that the combinatorial
definition of a ribbon operator can be generalized and used to give a
formula analogous to Theorem \ref{ribop}.

Let $V \in Hom(\Lambda,\Lambda)$ be an operator that does not involve
the parameter $q$ and let $P \in \Lambda$  also not include the
parameter $q$.  By setting $q=0$ in the expression $\htq{V} P[X]$, 
\begin{equation}\htq{V} P[X] \coeff_{q=0} = V^Y P[q X +(1-q)Y] \coeff_{Y=X} \coeff_{q=0}=
V P[X].\end{equation}
We remark that the highest power of $q$ that appears in this expression is
the degree of the symmetric function $P$.  By acting $\htq{V}$ on a 
Schur function, it can be seen that
\begin{equation}\htq{V} s_\la[X] = \sum_{\mu \subseteq \la} q^{|\mu|} V^Y ( s_\mu[X-Y]
s_{\la\slash\mu}[Y] ) \coeff_{Y=X}.\end{equation}
The coefficient of $q^{|\la|}$ in this expression will be the term
\begin{equation}\htq{V} s_\la[X] \coeff_{q^{|\la|}} = V^Y s_{\la}[X-Y] \coeff_{Y=X} 
= \ahat{V} s_\la[X].\end{equation}
The coefficients of $q^k$ may be interpreted then as a discrete interpolation
between $V$ and $\ahat{V}$.  Define notation for the coefficient of
$q^k$ of this expression so that
\begin{equation}\htq{V} P[X] \coeff_{q^k} = \ahat{V}^{(k)} P[X] =
\sum_{\la \vdash k} V^Y ( s_\la[X-Y] (s_\la^\perp P)[Y]) \coeff_{Y=X}. \label{khat}\end{equation}

By linearity, this notation may be extended to any operator $V$ that may
now depend on the parameter $q$.  This yields the following proposition.

\begin{prop} \label{khatprop} For $V \in Hom( \Lambda, \Lambda)$,
\begin{equation}\htq{V} = \sum_{k \geq 0} q^k \ahat{V}^{(k)},\end{equation}
where $\ahat{V}^{(k)}$ is defined in equation (\ref{khat}).
\end{prop}

$\ahat{S^R}^{(k)}$ can be developed in detail
thereby giving a combinatorial method for calculating the coefficient
of $q^i$ in a Macdonald polynomial.

Define a notion of a generalized ribbon operator that starts
with a ribbon $R \models m$ with $R = \la \slash \la^{rc}$ and associate
with this a sequence $v = (v_1, v_2,
\ldots, v_m)$ with $v_i \geq 0$. Generalize the notion of a ribbon
by setting the `thickness' of the $i^{th}$ cell of the ribbon to be
$v_i+1$ so that when the sequence consists of $m$ zeros this gives the
standard ribbon.

Let $\ell= \ell(\la)$ and say
that $D(R) = \{ i_1 > i_2 > \cdots > i_{\ell-1}\}$ and 
$\{1, \ldots, m\} - D(R) = \{j_1 < j_2 < \cdots < j_{\la_1}\}$.
Let $\alpha = \la^{rc} - (v_{i_{1}+1}, v_{i_{2}+1}, \ldots,
v_{i_{\ell-1}+1})$ (as vectors) and $\beta' = \la' + (v_{1}, v_{j_1+1},
v_{j_2+1}, \ldots, v_{j_{\la_1-1}+1})$ (neither $\alpha$ nor $\beta'$ are necessarily
partitions). Then set $S^{(R,v)} = (-1)^{|\la^{rc}|-|\alpha|}
s_\alpha^\perp \Sw_{\beta_1'} \Sw_{\beta_2'}
\cdots \Sw_{\beta_{\la_1}'}$.  Call $S^{(R,v)}$ a generalized ribbon
operator.

The formulation of these operators leads to a simple construction
with a picture: draw the original ribbon and place $v_i$ cells either
to the left of the $i^{th}$ cell if $i-1$ is a descent of the ribbon or above 
the cell if it is not. 
$\alpha$ is the sequence representing the space underneath the diagram and
$\beta'$ is the sequence representing the heights of the columns of the
diagram.  The sign represents the number of cells that are `underneath' the
ribbon. We present a couple of examples to give a better picture of these
truly combinatorial constructions.

\begin{ex} 
\end{ex}  
Consider the ribbon $R= \young{&\cr\blk&\cr\blk&\cr}$ ~of size $4$
with $D(R) = \{3,2\}$.  If $v=(0,0,0,0)$, then $S^{(R,v)} = S^{R}$. 
 If $v$ is one of $(1,0,0,0)$, $(0,1,0,0)$,
$(0,0,1,0)$, $(0,0,0,1)$  then $S^{(R,v)}$ is equivalent to the
following ribbon operators (respectively)

$$\young{\cr.&.\cr\blk&.\cr\blk&.\cr} = s_{11}^\perp \Sw_4 \Sw_3 \hskip
.3in 
\young{\blk&\cr.&.\cr\blk&.\cr\blk&.\cr} = s_{11}^\perp \Sw_3 \Sw_4 \hskip
.3in  -\young{.&.\cr&.\cr\blk&.\cr} = -s_{1}^\perp \Sw_3 \Sw_3 \hskip .3in 
-\young{.&.\cr\blk&.\cr&.\cr} = -s_{01}^\perp \Sw_3 \Sw_3 $$
The second and fourth generalized ribbons are $0$.  The second because it
contains the operation of adding a column of size $3$ on a column of size $4$,
and the fourth because a row of size $0$
is added on a row of size $1$ in the skew
part of the operator.

\begin{ex} 
\end{ex}  
 Let $R = \young{&&\cr\blk&\blk&\cr\blk&\blk&&&\cr\blk&\blk&\blk&\blk&\cr
\blk&\blk&\blk&\blk&\cr}$.  Now let $v = (1,1,0,1,0,2,0,1,0)$
This is represented by the following picture where a dot is placed in
each of the cells representing the original ribbon and there are $v_i$
cells either to the left of the $i^{th}$ cell if $i-1 \in
D(R)$ or above if $i-1 \notin D(R)$ (and the $v_1$ cells always go above the first cell in the ribbon).

$$\young{&\cr.&.&.&\cr\blk&&.&\cr
\blk&\blk&.&.&.\cr\blk&\blk&\blk&&.\cr
\blk&\blk&\blk&\blk&.\cr} = s_{4321}^\perp \Sw_6 \Sw_6 \Sw_5 \Sw_5 \Sw_3$$

Representing these operators with a diagram of this sort works fine if
$i-1$ is a descent and $v_i$ is so large that it creates a negative index
in $\alpha$.  Interpret this to mean that skewing by a Schur
function with a negative index kills the term and the result is $0$.

Note also that some `straightening' using the relation $S_m S_n = - S_{n-1} S_{m+1}$ may be
necessary.

\begin{ex} 
\end{ex}  
$R$ is as above,
but $v = (1,1,0,1,0,4,0,1,5)$.  Then $R$ can be represented by the image
$$-\young{\blk&\blk&\blk&\blk&\cr
\blk&&&\blk&\cr\blk&.&.&.&\cr\blk&\blk&&.&\cr\blk&\blk&\blk&.&.&.\cr
\blk&\blk&\blk&\blk&&.\cr
& & & & &.\cr} =- s_{-1,3,2,1}^\perp \Sw_6 \Sw_6 \Sw_5 \Sw_7 \Sw_3 =
 s_{-1,3,2,1}^\perp \Sw_6 \Sw_6 \Sw_6 \Sw_6 \Sw_3 =
s_{2,1,1,1}^\perp \Sw_6 \Sw_6 \Sw_6 \Sw_6 \Sw_3
= \young{&&&\cr&&&\cr\blk&&&\cr\blk& &&&\cr\blk& & &&\cr
\blk&\blk& & &\cr}$$
The final image comes from first straightening the columns of the
generalized ribbon and then straightening
Schur function that one skews by with appropriate sign changes.

Generalized ribbon operators are related to the original notion of
a ribbon operator by the following easily statable theorem.

\begin{thm} \label{genribrule}
Let $R$ be a ribbon of size $m$  and $k \geq 0$ an integer.
\begin{equation}\ahat{S^R}^{(k)} = \sum_{v} S^{(R,v)} e_v^\perp,\label{khatrib}\end{equation}
where  the sum is over all sequences $v$ having length $m$ and whose
sum is $k$ and the condition that $v_i \geq 0$ and $e_v$ is the elementary
symmetric function indexed by the sequence $v$.  Let
\begin{equation}H_{1^m}^{qt} = \sum_{R \models m} \sum_{v} q^{|v|} t^{comaj(R)} S^{(R,v)} 
e_v^\perp, \label{genrib}\end{equation}
where here the sum is over all sequences $v$ having length $m$ and non-negative
entries.
Then  $H_{1^m}^{qt} H_\mu[X;q,t] = H_{1^m|\mu}[X;q,t]$ for
$\ell(\mu) \leq m$.
\end{thm}

This theorem is a
combinatorial rule for computing Macdonald symmetric functions.
Before we present the proof, we give an example of how this theorem
works.

\begin{ex} Computation of a Macdonald symmetric function with generalized ribbons
\end{ex}  
 We will use formula (\ref{genrib}) to compute $H_{222}[X;q,t]$.  This
is a long and involved example, but it demonstrates the power of this
 this recurrence since with a reasonable amount of
work one can calculate a Macdonald polynomial of size $6$ or higher by hand.

Start with the formula for
$H_{111}[X;q,t] =
\young{\cr\cr\cr}+(t+t^2)
\young{\cr&\cr} +t^3 \young{&&\cr}$ (this may be calculated 
by acting $H_{1^3}^{qt}$ on $1$).

The sum in equation (\ref{genrib}) over $v$ is finite because only
terms such that
$|v|$ is less than or equal to the degree of the symmetric function that
is being acted on are needed.  In this computation quite a few operators are necessary.

We list all of the relevant operators (those which are non-zero) and place
a dot in the cells that consist of a the core of the operator so that it is
easy to read the sequence $v$ from the picture.  The sign associated to each
picture of the operator is $-1$ to the power of the number of cells under 
the ribbon.
$$\young{.\cr.\cr.\cr} 
\hskip .2in \young{\cr.\cr.\cr.\cr} 
\hskip .2in \young{\cr\cr.\cr.\cr.\cr}
\hskip .2in \young{\cr\cr\cr.\cr.\cr.\cr}
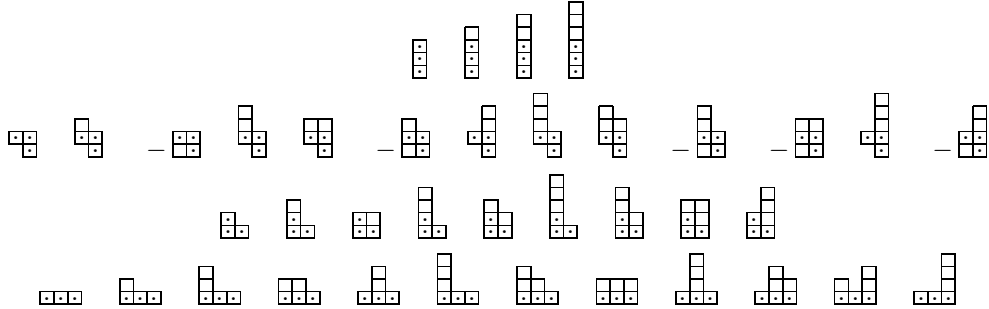
$$
$$\young{.&.\cr\blk &.\cr}
\hskip .2in \young{\cr.&.\cr\blk &.\cr} 
\hskip .2in -\young{.&.\cr &.\cr}
\hskip .2in \young{\cr\cr.&.\cr\blk &.\cr}
\hskip .2in \young{&\cr.&.\cr\blk &.\cr}
\hskip .2in -\young{\cr.&.\cr &.\cr}
\hskip .2in \young{\blk&\cr\blk&\cr.&.\cr\blk &.\cr}
\hskip .2in \young{\cr\cr\cr.&.\cr\blk &.\cr}
\hskip .2in \young{\cr&\cr.&.\cr\blk &.\cr}
\hskip .2in -\young{\cr\cr.&.\cr &.\cr}
\hskip .2in -\young{&\cr.&.\cr &.\cr}
\hskip .2in \young{\blk&\cr\blk&\cr\blk&\cr.&.\cr\blk &.\cr}
\hskip .2in -\young{\blk&\cr\blk&\cr.&.\cr &.\cr}$$
$$\young{.\cr. &.\cr} 
\hskip .2in \young{\cr.\cr. &.\cr} 
\hskip .2in \young{.&\cr. &.\cr} 
\hskip .2in \young{\cr\cr.\cr. &.\cr}
\hskip .2in \young{\cr.&\cr. &.\cr}
\hskip .2in \young{\cr\cr\cr.\cr. &.\cr}
\hskip .2in \young{\cr\cr.&\cr. &.\cr}
\hskip .2in \young{&\cr.&\cr. &.\cr}
\hskip .2in \young{\blk&\cr\blk&\cr.&\cr. &.\cr}$$
$$\young{.&.&.\cr} 
\hskip .2in \young{\cr.&.&.\cr}  
\hskip .2in \young{\cr\cr.&.&.\cr} 
\hskip .2in \young{&\cr.&.&.\cr} 
\hskip .2in \young{\blk&\cr\blk&\cr.&.&.\cr}
\hskip .2in \young{\cr\cr\cr.&.&.\cr}
\hskip .2in \young{\cr&\cr.&.&.\cr}
\hskip .2in \young{&&\cr.&.&.\cr}
\hskip .2in \young{\blk&\cr\blk&\cr\blk&\cr.&.&.\cr}
\hskip .2in \young{\blk&\cr\blk&&\cr.&.&.\cr} 
\hskip .2in \young{\blk&\blk&\cr&\blk&\cr.&.&.\cr}
\hskip .2in \young{\blk&\blk&\cr\blk&\blk&\cr\blk&\blk&\cr.&.&.\cr}$$

To complete this computation, calculate $e_\la^\perp$ on the symmetric function
$H_{111}[X;t]$ for $|\la|\leq 3$.  This is given by the following list

\begin{equation}e_1^\perp H_{111}[X;t] = 
\left(1+t+{t}^{2}\right ) \left( t s_{{2}}+s_{{1,1}}\right)\end{equation}
\begin{equation}e_2^\perp H_{111}[X;t] = \left (1+t+{t}^{2}\right )s_{{1}}\end{equation}
\begin{equation}e_{11}^\perp H_{111}[X;t] = \left (1+2\,t+2\,{t}^{2}+{t}^{3}\right
)s_{{1}}\end{equation} 
\begin{equation}e_{3}^\perp H_{111}[X;t] = 1
\end{equation}
\begin{equation}e_{21}^\perp H_{111}[X;t] = 1+t+{t}^{2}
\end{equation}
\begin{equation}e_{111}^\perp H_{111}[X;t] = 1+2\,t+2\,{t}^{2}+{t}^{3}
\end{equation}

The computation proceeds as follows.  The coefficient of $q^0$ is just the
Hall-Littlewood symmetric function $H_{222}[X;t]$, calculated by acting
$\young{.\cr.\cr.\cr} + t \young{.&.\cr\blk &.\cr} + 
t^2 \young{.\cr. &.\cr} + t^3 \young{.&.&.\cr}$
on $H_{111}[X;t]$.
So the coefficient of $q^0$ is
\begin{equation} s_{222} + (t +t^2)s_{321}
+ t^3 s_{33} +
t^3 s_{411} + (t^2+t^3+t^4) s_{42}  +  (t^4+t^5 )s_{51} + t^6 s_{6}\end{equation}

The coefficient of $q^1$ is the operator $\young{\cr.\cr.\cr.\cr}
+ t \young{\cr.&.\cr\blk &.\cr}-t \young{.&.\cr &.\cr}+t^2 \young{\cr.\cr. &.\cr} 
+ t^2\young{.&\cr. &.\cr}+t^3\young{\cr.&.&.\cr}$ when it acts on the symmetric
function $\left(1+t+{t}^{2}\right ) \left( t s_{{2}}+s_{{1,1}}\right).$
\begin{equation}\left(1+t+{t}^{2}\right ) (t^4 s_{51}+ (t^2+t^3) s_{411} + t^3 s_{42}
+ (t+t^2) s_{321}+ t^2 s_{33} + ts_{3111} + s_{2211}    )\end{equation}

The coefficient of $q^2$ comes from two components, $\young{\cr\cr.\cr.\cr.\cr}+
t \young{\cr\cr.&.\cr\blk &.\cr}+t \young{\blk&\cr\blk&\cr.&.\cr\blk &.\cr}+
t^2 \young{\cr\cr.\cr. &.\cr}+t^3 \young{\cr\cr.&.&.\cr}+t^3\young{\blk&\cr\blk&\cr.&.&.\cr}$
when it acts on the symmetric function $e_2^\perp H_{111}[X;q,t]$, and
$t \young{&\cr.&.\cr\blk &.\cr}
 -t \young{\cr.&.\cr &.\cr} + t^2 \young{\cr.&\cr. &.\cr} +
t^3 \young{&\cr.&.&.\cr}$ when it acts on the symmetric function 
$e_{11}^\perp H_{111}[X;q,t]$.
The first part is
\begin{equation}\left ({t}^{2}+t+1\right ) \left (
-{t}^{3}s_{{42}}+{t}^{3}s_{{411}}+\left (t+t^2\right )s_{{3111}}-
ts_{{222}}+s_{{2211}}t+s_{{21111}}
 \right) \end{equation}
and the second is
\begin{equation}\left (1+2\,t+2\,{t}^{2}+{t}^{3}\right ) \left( 
{t}^{3}s_{{42}}+{t}^{2}s_{{321}}+ts_{{222}}
\right)\end{equation}
The sum of these two quantities is
\begin{equation}\left ({t}^{2}+t+1\right )\left
({t}^{4}s_{{42}}+{t}^{3}s_{{411}}+\left ({t}^{2}+{t}^{3}\right )s_{{321
}}+\left (t+{t}^{2}\right )s_{{3111}}+{t}^{2}s_{{222}}+ts_{{2211}}+
s_{{21111}}
\right )
\end{equation}

The coefficient of $q^3$ comes from three different operators acting each
on a different constant.  The first operator is $\young{\cr\cr\cr.\cr.\cr.\cr}+
t \young{\cr\cr\cr.&.\cr\blk &.\cr} + t \young{\blk&\cr\blk&\cr\blk&\cr.&.\cr\blk &.\cr}
+t^2 \young{\cr\cr\cr.\cr. &.\cr}+t^2 \young{\blk&\cr\blk&\cr.&\cr. &.\cr}+
t^3 \young{\cr\cr\cr.&.&.\cr}+t^3 \young{\blk&\cr\blk&\cr\blk&\cr.&.&.\cr}+
t^3 \young{\blk&\blk&\cr\blk&\blk&\cr\blk&\blk&\cr.&.&.\cr}$ when it acts on  $1$,
the second is $t \young{\cr&\cr.&.\cr\blk &.\cr} - t \young{\cr\cr.&.\cr &.\cr} 
-t \young{\blk&\cr\blk&\cr.&.\cr &.\cr}+t^2 \young{\cr\cr.&\cr. &.\cr}
+ t^2 \young{&\cr.&\cr. &.\cr}+t^3 \young{\cr&\cr.&.&.\cr}
+t^3 \young{\blk&\cr\blk&&\cr.&.&.\cr}+t^3 \young{\blk&\blk&\cr&\blk&\cr.&.&.\cr}$ when it acts on $1+t+t^2$, and the
third is $-t \young{&\cr.&.\cr &.\cr} +t^3 \young{&&\cr.&.&.\cr}$ when it acts on
the symmetric function $2\,t+2\,{t}^{2}+{t}^{3}+1$.  These three parts are

\begin{equation}{t}^{3}s_{{33}}-{t}^{3}s_{{321}}+{t}^{3}s_{{3111}}-\left (t+{t}^{2}
\right )s_{{222}}+\left (t+{t}^{2}\right )s_{{21111}}+s_{{111111}}\end{equation}

\begin{equation}\left(1+t+t^2\right)\left(
-2\,{t}^{3}s_{{33}}+{t}^{3}s_{{321}}+\left ({t}^{2}+2 t\right )s_{{222}}
+{t}^{2}s_{{2211}} \right)\end{equation}

\begin{equation}\left (2\,t+2\,{t}^{2}+{t}^{3}+1 \right)\left( -t s_{222} + t^3
s_{33} \right)\end{equation} 
The sum of these three quantities is
\begin{equation}{t}^{6}s_{{33}}+\left ({t}^{4}+{t}^{5}\right
)s_{{321}}+{t}^{3}s_{{3111}}+{t}^{3}s_{{222}}+\left ({t}^{2}+t+1\right
){t}^{2}s_{{2211}} +\left (t+{t}^{2}\right )s_{{21111}}+s_{{111111}}
 \end{equation} which is the coefficient of $q^3$ in $H_{222}[X;q,t]$.

Clearly, an enormous amount of simplification occurs when arriving
at a final expression for $H_\mu[X;q,t]$.  An eventual goal of a
combinatorial recurrence on the $q,t$-Kostka coefficients will be to
arrive at a combinatorial interpretation for them in terms of standard
tableaux.  Even if this recurrence turns out to be too complicated, 
these techniques (in particular, Theorem \ref{Hqtcreation} and Theorem \ref{ribop}) 
can certainly be used to derive
many other recurrences for the coefficients.

Before presenting the proof, we will need a few lemmas that come from the
derivation of the ribbon operator.  We state them without
proof and refer the reader to [Za1].

\begin{lem}\label{lem1} For any operator $V$, 
$\ahat{\ahat{V}S_m} = \sum_{j\geq 0}(-1)^{m-j} h_j V \Sw_{m-j}$ \end{lem}

\begin{lem}\label{lem2}
 $\ahat{\ahat{s_{\la}^\perp}S_{-m}} = 
(s_{(m,\la)})^\perp$.  \end{lem}

\bproof{ (of Theorem \ref{genribrule})} If $a \geq 0$ and $v$ is a list, then we
denote $v$ with $a$ prepended
(resp. appended) by $(a,v)$ (resp. $(v,a)$).

Let $R^+$ be a ribbon of size $m+1$ that does not
have $m$ as a descent.  Also let $R$ be the ribbon of size $m$
such that $D(R) = D(R^+)$.  By the definition of $S^{(R,v)}$,
notice that $S^{(R^+,(v,a))} = S^{(R,v)} \Sw_{1+a}.$

Now let $R_+$ be a ribbon of size $m+1$
such that $m$ is a descent.  Let $R = \la\slash \la^{rc}$
be the ribbon of size $m$
such that $D(R) \cup \{m\} = D(R_+)$ then remark that
$R_+ = (\la_1, \la) \slash (\la_1-1, \la^{rc})$.
If $S^{(R,v)} = (-1)^{|v|} s_{\alpha}^\perp \Sw_{\beta'}$,
then $S^{(R_+, (v,a))} = (-1)^{|v|+a} 
s_{(\la_1-1-a,\alpha)}^\perp \Sw_{(\la_1,\beta)'}$.
It follows from Lemma \ref{lem1} and \ref{lem2} and the commutation
relation $\Sw_m \Sw_n = - \Sw_{n-1} \Sw_{m+1}$ that
\begin{eqnarray}S^{(R_+, (v,a))} &=& (-1)^{|\la^{rc}|-|\alpha|+a} 
\ahat{\ahat{s_{\alpha}^\perp}S_{1+a-\la_1}} \Sw_{(\la_1,\beta)'}\nonumber\\
&=& (-1)^{|\la^{rc}|-|\alpha|+a} \sum_{j\geq 0}(-1)^{1+a-\la_1-j} h_j s_{\alpha}^\perp \Sw_{1+a-\la_1-j}
\Sw_{(\la_1,\beta)'}\nonumber\\
&=& (-1)^{|\la^{rc}|-|\alpha|+a} \sum_{j\geq 0}(-1)^{1+a-j} h_j s_{\alpha}^\perp 
\Sw_{\beta'}\Sw_{1+a-j}\\
&=& (-1)^a \ahat{\ahat{S^{(R,v)}} S_{1+a}}\nonumber\end{eqnarray}

Use these two relations to give an inductive derivation of the following
plethystic form of the operator $S^{(R,v)}$.  By carrying out nearly the
exact same calculation (and using identical notation for
$ Z_{D(R)^c+1}$) as given in Proposition \ref{plethrib}, derive that
\begin{eqnarray}S^{(R,v)} P[X] &=& (-1)^{m- |D(R)|+|v|} P[ X+ Z^*] 
\Omega[-(z_1 + Z_{D(R)^c+1}) X]\nonumber\\
& &\prod_{1 \leq i < j \leq m}
(1-z_j/z_i) \coeff_{z_1^{1+v_1} z_2^{1+v_2} \cdots z_m^{1+v_m}}.\end{eqnarray}

Now consider a formula for $\htq{S^R}$.  Using the same calculation for Theorem
\ref{plethH1mqt} and the equation given in Proposition \ref{plethrib}, demonstrate that
\begin{equation}\htq{S^R} = (-1)^{m-|D(R)|} P[ X+ (1-q)Z^*] \Omega[-(z_1 + Z_{D(R)^c+1}) X] \prod_{1 \leq i < j \leq m}
(1-z_j/z_i) \coeff_{z_1 z_2 \cdots z_m}.\end{equation}

The coefficient of $q^k$ in this formula is 

\begin{eqnarray}\ahat{S^R}^{(k)} &=& (-1)^{m-|D(R)|} \sum_{\la\vdash k} f_\la[- Z^{*}] e_\la^\perp P[ X+ Z^*] \Omega[-(z_1 + Z_{D(R)^c+1}) X] \prod_{1 \leq i < j \leq m}
(1-z_j/z_i) \coeff_{z_1 z_2 \cdots z_m
}\nonumber\\
&=&(-1)^{m-|D(R)|+k} \sum_{\la\vdash k} m_\la[Z^{*}] e_\la^\perp P[ X+ Z^*] \Omega[-(z_1 + Z_{D(R)^c+1}) X] \prod_{1 \leq i < j \leq m}
(1-z_j/z_i) \coeff_{z_1 z_2 \cdots z_m
}.\nonumber\end{eqnarray}

By expanding $m_\la[Z^{*}]$ as $\sum_{v \sim \la} z^{-v}$ we see clearly that this
is equivalent to  equation  (\ref{khatrib}).  
The formula stated for the operator $H_{1^m}^{qt}$ follows
from this derivation, Theorem \ref{Hqtcreation} and Proposition \ref{khatprop}. \endofproof

\section{Bibliography}

\article BGHT|F.~Bergeron, A.~M.~Garsia, M.~Haiman and G.~P.~Tesler|Identities and positivity
conjectures for some remarkable operators in the theory of symmetric functions|
Methods of Analysis and Applications|6|(1999)|58 pp|

\article G|A.~M.~Garsia|Orthogonality of Milne's polynomials and
raising operators|Discrete Math.|99|(1992)|247--264|

\article GHT|A.~M.~Garsia, M.~Haiman and G.~P.~Tesler|Explicit 
Plethystic Formulas for the Macdonald q,t-Kostka Coefficients|The Andrews Festschrift (Maratea, 1998),
S\'em. Lothar. Combin.|42|(1999)|45 pp|

\article GZ|A.~M.~Garsia and M.~Zabrocki|Polynomiality of the 
$q,t$-Kostka Revisited|Algebraic combinatorics and computer
science, Springer Italia, Milan||(2001)|473--491|

\article J|N.~Jing|Vertex operators and Hall-Littlewood symmetric
functions|Adv. Math.|87|(1991)|226--248|

\unarticle K|A.N.~Kirillov|Ubiquity of Kostka polynomials|math.QA/9912094||

\article KN1|A.N.~Kirillov and M.~Noumi|$q$-difference raising operators for Macdonald polynomials and the integrality of transition coefficients. 
Algebraic methods and $q$-special functions| 
CRM Proc. Lecture Notes (Montr\'eal, QC, 1996)
|22|, Amer. Math. Soc., Providence, RI (1999)|227--243|

\article KN2|A.N.~Kirillov and M.~Noumi|Affine Hecke algebras and raising 
operators for Macdonald polynomials|
Duke Math. J|93|(1998), no. 1|1--39|

\unarticle LLM|L.~Lapointe, A.~Lascoux, J.~Morse|Tableau atoms and a 
new Macdonald positivity conjecture|preprint||

\unarticle LM1|L.~Lapointe, J.~Morse|Schur function identities, 
their t-analogs, and k-Schur irreducibility|preprint||

\unarticle LM1|L.~Lapointe, J.~Morse|Schur function analogs for 
a filtration of the symmetric function space|preprint||

\article LV1|L.~Lapointe and L.~Vinet|A short proof of the integrality of
the Macdonald $(q,t)$-Kostka coefficients|
Duke Math. J.|91|(1998), no. 1|205--214|

\article LV2|L.~Lapointe and L.~Vinet|Rodrigues formulas for the Macdonald 
polynomials|Adv. Math.|130|(1997), no. 2|261--279|

\article M1|I. G. Macdonald|A new class of symmetric functions|
Actes du $20^e$ S\'eminaire Lotharingien| 
|Publ. I.R.M.A. Strasbourg (1988)|
131-171|

\livre M2|I.~G.~Macdonald|Symmetric Functions and Hall
Polynomials, Oxford Mathematical Monographs|second edition, 
Oxford Univ. Press, 1995|

\article SW|M.~Shimozono and J.~Weyman|
Graded characters of modules supported in the
closure of a nilpotent conjugacy class|European J.
Combin.|21|(2000), no. 2|257--288|

\article SZ|M.~Shimozono and M.~Zabrocki|Hall-Littlewood vertex operators and
generalized Kostka polynomials|Adv. Math|158|(2001)|66--85|

\article Za1|M.~Zabrocki|Ribbon Operators and 
Hall-Littlewood Symmetric Functions|Adv. Math.|156|(2000), no. 1|33--43|

\article Za2|M.~Zabrocki|Vertex Operators for Standard Bases of the
Symmetric Functions|J. of Alg. Comb.|13|(2000), no. 1|83-101|

\livre Ze|A.~V.~Zelevinsky|Representations of finite classical groups:
a Hopf algebra approach|Springer Lecture Notes, 869, 1981|
\end{document}